\documentclass[11pt,reqno]{amsart}

\usepackage[T1]{fontenc}
\usepackage{lmodern}
\usepackage{amssymb}
\usepackage{mathtools}
\usepackage[a4paper,margin=2.7cm]{geometry}
\usepackage{enumitem}
\usepackage{xcolor}
\usepackage{xurl}
\usepackage{hyperref}
\hypersetup{
  colorlinks=true,
  linkcolor=blue!55!black,
  citecolor=blue!55!black,
  urlcolor=blue!55!black,
  pdftitle={Measure-class classification of achievable Cantorvals},
  pdfauthor={HanYu Dong}
}

\newtheorem{theorem}{Theorem}[section]
\newtheorem{lemma}[theorem]{Lemma}
\newtheorem{proposition}[theorem]{Proposition}
\newtheorem{corollary}[theorem]{Corollary}
\newtheorem*{maintheorem}{Main Theorem}
\theoremstyle{definition}
\newtheorem{definition}[theorem]{Definition}
\newtheorem{example}[theorem]{Example}
\theoremstyle{remark}
\newtheorem*{remark*}{Remark}

\DeclareMathOperator{\supp}{supp}
\DeclareMathOperator{\intr}{int}
\DeclareMathOperator{\cl}{cl}
\DeclareMathOperator{\Homeo}{Homeo}
\DeclareMathOperator{\MC}{MC}
\DeclareMathOperator{\dist}{dist}
\DeclareMathOperator{\diam}{diam}
\DeclareMathOperator{\Spec}{Spec}
\DeclareMathOperator{\Hom}{Hom}

\newcommand{\R}{\mathbb{R}}
\newcommand{\mc}{\approx_{\mathrm{mc}}}

\begin{document}

\title[Measure classes of achievable Cantorvals]{Measure-class classification of achievable Cantorvals}

\author{HanYu Dong}
\address{School of Mathematics, Southwest Minzu University, Chengdu, Sichuan, China}
\email{shudaiziyicc@qq.com}

\date{November 19, 2026}



\begin{abstract}
An achievable Cantorval is a Cantorval that is the set of subsums of a convergent series of positive terms. We construct an explicit achievable Cantorval whose boundary has Lebesgue measure at least $2\prod_{j\ge 1}(1-2^{-j})=0.5775\ldots$. This answers a question of G{\l}{\k{a}}b and Prus-Wi{\'s}niowski and gives a counterexample to the conjecture of Nowakowski and Prus-Wi{\'s}niowski that every achievable Cantorval is standard. The example is generated by a hybrid series of Pratsiovytyi, Vynnyshyn and Vynnyshyn. The lower bound comes from an isolated-digit criterion, which produces a Cantor subset of the boundary whose measure is an explicit infinite product; after rescaling, the boundary measure can be made arbitrarily close to $1$. We also prove that the boundary of every achievement set is either Lebesgue null or has positive measure near each of its points. Consequently, up to homeomorphisms that preserve Lebesgue null sets in both directions, there are exactly two achievable Cantorvals, distinguished by whether the boundary is null. For boundaries of positive measure, the homeomorphism can be chosen so that it changes the normalized boundary measure of every Borel set by a factor arbitrarily close to $1$. For general Cantorvals we classify the case of a relatively clopen boundary support and show that further invariants appear beyond this case. Finally, the null-set preserving self-homeomorphisms of a Cantorval always form a meager subgroup of its homeomorphism group. This subgroup is dense exactly when the support of Lebesgue measure on the boundary is empty or the whole boundary; in particular, it is dense for every achievable Cantorval.
\end{abstract}

\maketitle

\section{Introduction}

For a summable sequence $x=(x_n)$ of positive numbers, its \emph{achievement set} (or \emph{set of subsums}) is
\[
E(x)=\Bigl\{\sum_{n=1}^{\infty}\varepsilon_n x_n : \varepsilon_n\in\{0,1\}\Bigr\}.
\]
Achievement sets are compact. By the theorem of Guthrie and Nymann~\cite{GN} (see also Nymann and S{\'a}enz~\cite{NS}), each of them is a finite union of closed intervals, a Cantor set, or a Cantorval. In this paper a \emph{Cantorval} is a nonempty compact set $X\subset\R$ that is the closure of its interior and whose boundary is a Cantor set. By \cite[Theorem~2.4]{MNP} this agrees with the usual definition, which goes back to Mendes and Oliveira~\cite{MO}. All Cantorvals are homeomorphic to each other \cite[Appendix]{MO}; order-theoretic descriptions of Cantorvals and of other compact subsets of the line can be found in \cite{MO,BKP,KP}. For background on achievement sets we refer to the surveys \cite{BFP} and \cite{GPW}.

Since all Cantorvals are homeomorphic, finer invariants are needed to tell them apart, and Lebesgue measure provides natural ones. In particular, the Lebesgue measure of the boundary of an achievable Cantorval has been the subject of several recent problems. Nowakowski and Prus-Wi{\'s}niowski~\cite{NPW} proved that standard Cantorvals (see Section~\ref{ss:standard}) have null boundaries, and conjectured that every achievable Cantorval is standard \cite[Conjecture~11]{NPW}. In the preprint version of their survey, G{\l}{\k{a}}b and Prus-Wi{\'s}niowski \cite[Problems~21 and~22]{GPW} asked whether an achievable Cantorval can have a boundary of positive Lebesgue measure, and whether all achievable Cantorvals are standard. The survey also discusses the relations between tail inequalities, interval components and the geometry of the boundary; see also~\cite{PWP}. The boundary measure and the dimension of the Guthrie--Nymann Cantorval and of its generalizations are studied in~\cite{PK}.

We study these questions through homeomorphisms that preserve Lebesgue null sets in both directions; we call them \emph{measure-class homeomorphisms} (Definition~\ref{def:mc}). For a Cantorval $X$ we write
\[
S(X)=\supp(\lambda|_{\partial X}),
\]
where $\lambda$ denotes Lebesgue measure. This set is invariant under measure-class homeomorphisms (Lemma~\ref{lem:supp}). Our main result is the following.

\begin{maintheorem}
The following assertions hold.
\begin{enumerate}[label=\textup{(\roman*)},leftmargin=*]
\item For the achievement set $E$ of every summable sequence of positive numbers,
\[
\lambda(\partial E)>0 \implies \supp(\lambda|_{\partial E})=\partial E.
\]
\item There is an explicit achievable Cantorval $E_*$ such that
\[
\lambda(\partial E_*)\ge 2\prod_{j=1}^{\infty}(1-2^{-j})=0.5775761901732048\ldots.
\]
Moreover, there are achievable Cantorvals with convex hull $[0,1]$ whose boundaries have measure arbitrarily close to $1$.
\item The achievable Cantorvals form exactly two measure-class homeomorphism classes, distinguished by whether the boundary is null or has positive measure. Any two members of the same class are related by a measure-class homeomorphism that is the restriction of an increasing homeomorphism between their convex hulls, affine on every interval component and on every gap. In the positive class, this homeomorphism can moreover be chosen to distort the normalized boundary measures by a factor of at most $L$, for any prescribed $L>1$.
\end{enumerate}
\end{maintheorem}

Parts (i), (ii) and (iii) are Theorem~\ref{thm:dichotomy}, Theorem~\ref{thm:hybrid} with Corollary~\ref{cor:close-to-one}, and Theorem~\ref{thm:two-classes}, respectively.

Part~(i) rests on the classical description of an achievement set $E$ as the intersection of the finite unions of intervals $H_n=\bigcup_{f\in F_n}[f,f+r_n]$, where $F_n$ is the set of subsums of the first $n$ terms and $r_n$ is the $n$-th remainder. Every gap of $E$ is already a gap of some $H_n$, and near each endpoint of such a gap the set $E$ coincides with a translate of the tail set $E_n$, the achievement set of $(x_k)_{k>n}$. We show that if $\partial E$ has positive measure, then so does the boundary of every tail set, and that this boundary measure accumulates at both endpoints of each tail set. Since the endpoints of gaps are dense in $\partial E$, every point of $\partial E$ belongs to the support.

For part~(ii) we use the hybrid series of Pratsiovytyi, Vynnyshyn and Vynnyshyn~\cite{PVV}. They are built from blocks of the form $Q_j^{-1}(b_j-1,\,b_j/2,\,b_j/4,\dots,2)$, where $b_j$ is a power of $2$ (see~\eqref{eq:blocks}), and their achievement sets are Cantorvals \cite[Theorem~4.2]{PVV}. Multiplied by $Q_j$, the subsums of the $j$-th block form a complete system of residues modulo $b_j$. These subsums contain exactly one pair of consecutive integers, namely $b_j-2$ and $b_j-1$, and any two other subsums differ by at least $2$. Discarding this pair at every scale leaves a Cantor set whose points are one-point components of the achievement set, and hence boundary points (Theorem~\ref{thm:isolated}). Its Lebesgue measure is an explicit infinite product, which is positive when $\sum_j 1/b_j<\infty$ (Theorem~\ref{thm:hybrid}). The Cantorval $E_*$ in~(ii) is the achievement set of the series of \cite[Example~2.5]{PVV}. A related construction of an achievable Cantorval whose boundary has full Hausdorff dimension is described in~\cite{PRISM}.

For part~(iii), recall that the bounded components of $\R\setminus\partial X$ are of two kinds: the interiors of the interval components of $X$, and the gaps of $X$. We regard them as cells of two colors. An increasing homeomorphism that maps the cells of $\partial X$ onto the cells of $\partial Y$ and preserves the colors maps $X$ onto $Y$ (Lemma~\ref{lem:lift}). If both boundaries are null, every such map is a measure-class homeomorphism. If both boundaries have positive measure, then by part~(i) the boundary measures have full support, and we pass to their distribution functions. This reduces the problem to the construction of an increasing bi-Lipschitz homeomorphism of $[0,1]$, with constants arbitrarily close to $1$, that maps two given disjoint countable dense sets onto two other such sets (Lemma~\ref{lem:interpolation}). This is a colored version of Akin's theorem on order homeomorphisms of Cantor spaces with arbitrarily small multiplicative distortion of measures \cite[Theorem~2.12]{Akin}, and places the classification in the setting of topological measure transport developed in \cite{OU,Akin}.

For an arbitrary Cantorval $X$, the support $S(X)$ may be a nonempty proper subset of $\partial X$. When $S(X)$ is relatively open in $\partial X$, the Cantorval splits into a part whose boundary has full support, a part with null boundary, and finitely many interval components joining the two parts. This leads to the classes $\mathsf{N}$, $\mathsf{F}$ and $\mathsf{M}_m$ ($m\ge 0$) of Theorem~\ref{thm:clopen}. When $S(X)$ is not open, further invariants appear (Proposition~\ref{prop:interface}). We also show that the measure-class homeomorphisms of a Cantorval always form a meager subgroup of its homeomorphism group, and that this subgroup is dense exactly when $S(X)$ is empty or equal to $\partial X$ (Theorems~\ref{thm:meager} and~\ref{thm:density}). For the classical relations between measure and category we refer to~\cite{Oxtoby}.

The paper is organized as follows. Section~\ref{sec:mc} introduces measure-class homeomorphisms and proves the transport theorems for Cantorvals whose boundaries are null or everywhere fat. Section~\ref{sec:support} proves the support dichotomy for achievement sets, and Section~\ref{sec:positive} constructs achievable Cantorvals with boundaries of positive measure. Section~\ref{sec:two-classes} derives the two-class theorem and the counterexample to the standardness conjecture. Sections~\ref{sec:clopen} and~\ref{sec:groups} treat mixed supports and homeomorphism groups. Appendix~\ref{app:cosmall} discusses maps defined on co-small subspaces, and Appendix~\ref{app:generic} identifies the analogous generic chain-homotopy relation for rational maps over a function field.

\section{Measure classes and homeomorphisms of Cantorvals}\label{sec:mc}

\begin{definition}\label{def:mc}
A \emph{measure-class space} is a triple $(X,\mathcal{T},\mathcal{I})$, where $\mathcal{T}$ is a topology on $X$ and $\mathcal{I}$ is a proper $\sigma$-ideal of subsets of $X$. A \emph{measure-class homeomorphism} from $(X,\mathcal{T},\mathcal{I})$ to $(Y,\mathcal{S},\mathcal{J})$ is a homeomorphism $h\colon X\to Y$ such that
\[
A\in\mathcal{I}\iff h(A)\in\mathcal{J}\qquad(A\subset X).
\]
We write $X\mc Y$ if such a homeomorphism exists. A compact set $X\subset\R$ of positive Lebesgue measure is always regarded as a measure-class space with the subspace topology and the ideal $\mathcal{N}(X)$ of subsets of $X$ of Lebesgue outer measure zero; this ideal is proper because $\lambda(X)>0$.
\end{definition}

Identities, inverses and compositions of measure-class homeomorphisms are measure-class homeomorphisms. If the ideals are the null ideals of finite Borel measures $\mu$ and $\nu$ on compact metric spaces, then a homeomorphism $h$ is a measure-class homeomorphism if and only if $h_*\mu$ and $\nu$ are mutually absolutely continuous. Topological spaces equipped with ideals are also studied in~\cite{JH}.

\begin{definition}\label{def:support}
For $B\subset X$, the \emph{support of $B$ relative to $\mathcal{I}$} is
\[
\supp_{\mathcal{I}}B=\{x\in X: U\cap B\notin\mathcal{I}\ \text{for every open neighborhood $U$ of $x$}\}.
\]
\end{definition}

The set $\supp_{\mathcal{I}}B$ is closed. If $X$ is second countable and $\mathcal{I}$ is the null ideal of a Borel measure $\mu$, then $\supp_{\mathcal{I}}B=\supp(\mu|_B)$ for every Borel set $B$. Moreover, $B\setminus\supp(\mu|_B)$ is $\mu$-null, since it is covered by the countably many basic open sets whose intersections with $B$ are null.

\begin{lemma}\label{lem:supp}
If $h\colon(X,\mathcal{T},\mathcal{I})\to(Y,\mathcal{S},\mathcal{J})$ is a measure-class homeomorphism, then
\[
h(\supp_{\mathcal{I}}B)=\supp_{\mathcal{J}}h(B)\qquad(B\subset X).
\]
\end{lemma}

\begin{proof}
The map $U\mapsto h(U)$ is a bijection from the open neighborhoods of $x$ onto those of $h(x)$, and $h(U\cap B)=h(U)\cap h(B)$ belongs to $\mathcal{J}$ if and only if $U\cap B$ belongs to $\mathcal{I}$.
\end{proof}

\subsection{Cantorvals, gaps and cells}

\begin{definition}\label{def:cantorval}
A \emph{Cantorval} is a nonempty compact set $X\subset\R$ such that $X=\cl(\intr X)$ and $\partial X$ is a Cantor set. A nondegenerate connected component of $X$ is called an \emph{interval component} of $X$, and a bounded connected component of $\R\setminus X$ is called a \emph{gap} of $X$. Put
\[
S(X)=\supp(\lambda|_{\partial X}).
\]
The boundary of $X$ is \emph{everywhere fat} if $S(X)=\partial X$, that is, if $\lambda(\partial X\cap U)>0$ for every open set $U$ meeting $\partial X$. More generally, a closed set $Z\subset\R$ is everywhere fat if $\supp(\lambda|_Z)=Z$.
\end{definition}

By \cite[Theorem~2.4]{MNP}, this definition of a Cantorval agrees with the usual one; see \cite{MO,BFP,GPW}. The two unbounded components of $\R\setminus X$ are sometimes called external gaps; in this paper gaps are always bounded unless stated otherwise.

For a Cantor set $Z\subset\R$, we call the bounded components of $\R\setminus Z$ the \emph{cells} of $Z$ and order them from left to right. (They are usually called the gaps of $Z$; we use a different word to avoid confusion with the gaps of a Cantorval.) If $Z=\partial X$ for a Cantorval $X$, then every cell of $Z$ is either the interior of an interval component of $X$ or a gap of $X$. Accordingly, we say that the cell has the \emph{interval color} or the \emph{gap color}.

\begin{lemma}[{see \cite{MO}}]\label{lem:homeo-invariance}
Let $h\colon A\to B$ be a homeomorphism between subsets of $\R$. Then $h(\intr A)=\intr B$. If $A$ and $B$ are compact, then $h$ maps $\partial A$ onto $\partial B$, and maps the interval components of $A$ and their endpoints onto the interval components of $B$ and their endpoints. In particular, a compact set homeomorphic to a Cantorval is a Cantorval.
\end{lemma}

These facts are standard, and we omit the proof. The first assertion is invariance of domain in dimension one; the others hold because homeomorphisms preserve connected components, and the endpoints of a nondegenerate compact interval are exactly the points whose removal leaves it connected.

\begin{lemma}\label{lem:cells}
The cells of a Cantor set $Z\subset\R$ form a countable dense linear order without endpoints. Their closures are pairwise disjoint and lie strictly between $\min Z$ and $\max Z$. If $Z=\partial X$ for a Cantorval $X$, then both colors are dense in this order, and every open interval meeting $\partial X$ contains the closures of cells of both colors.
\end{lemma}

\begin{proof}
A family of disjoint open intervals is countable. If two cell closures met, or if a cell closure contained $\min Z$ or $\max Z$, some point of $Z$ would be isolated. Between two cells there are two distinct points of $Z$, and since $Z$ is nowhere dense, there is a cell between them. Since $Z$ is perfect, there are cells arbitrarily close to $\min Z$ and to $\max Z$, so the order has no endpoints.

Now let $Z=\partial X$, and let $u<w<v$ be points of $\partial X$ in a given open interval. Since $w\in\partial X=\cl(\intr X)\setminus\intr X$, points of $\intr X$ and points of $\R\setminus X$ occur arbitrarily close to $w$, in particular in $(u,v)$; the cells containing them have their closures in $[u,v]$. Applying this argument between any two cells shows that both colors are dense.
\end{proof}

\begin{lemma}\label{lem:lift}
Let $Z,Z'\subset\R$ be Cantor sets and let $\sigma$ be an order isomorphism from the cells of $Z$ onto the cells of $Z'$. There is a unique increasing homeomorphism
\[
H\colon[\min Z,\max Z]\longrightarrow[\min Z',\max Z']
\]
that maps each cell $c$ affinely onto $\sigma(c)$. It satisfies $H(Z)=Z'$. If $Z=\partial X$ and $Z'=\partial Y$ for Cantorvals $X$ and $Y$, and $\sigma$ preserves the colors, then $H(X)=Y$.
\end{lemma}

\begin{proof}
Define $H$ on the closure of each cell $c$ as the increasing affine map onto the closure of $\sigma(c)$. By Lemma~\ref{lem:cells} these closures are pairwise disjoint, and since $Z$ and $Z'$ are nowhere dense, the unions of the cell closures are dense in the respective intervals. As $\sigma$ preserves the order, this defines an increasing bijection between two dense subsets of the two intervals, which also matches the endpoints of the intervals. Such a bijection extends uniquely to an increasing homeomorphism between the intervals. This homeomorphism maps the union of the cells of $Z$ onto the union of the cells of $Z'$, and hence maps $Z$ onto $Z'$. In the colored case it maps $\intr X$ onto $\intr Y$, and hence $X=\cl(\intr X)$ onto $Y$.
\end{proof}

Any two countable dense linear orders without endpoints, each partitioned into two dense subsets (colors), are isomorphic by a color-preserving order isomorphism. This follows from Cantor's back-and-forth argument: enumerate both orders and alternately extend a finite color-preserving partial isomorphism so that it contains the first element not yet used on either side; this is always possible because between any two elements, and below and above any element, there are infinitely many elements of each color. Together with Lemmas~\ref{lem:cells} and~\ref{lem:lift}, this gives an increasing homeomorphism between the convex hulls of any two Cantorvals that maps one onto the other. In particular, it recovers the theorem of Mendes and Oliveira that all Cantorvals are homeomorphic \cite[Appendix]{MO}.

\begin{lemma}\label{lem:pieces}
A homeomorphism $h\colon X\to Y$ between Cantorvals is a measure-class homeomorphism if and only if its restriction to $\partial X$ and its restriction to every interval component of $X$ preserve null sets in both directions. In particular, if $h$ is affine on every interval component, only the condition on the boundary has to be checked.
\end{lemma}

\begin{proof}
By Lemma~\ref{lem:homeo-invariance}, $h$ maps $\partial X$ onto $\partial Y$ and the interval components of $X$ onto those of $Y$. There are only countably many interval components, and every set $A\subset X$ is the union of $A\cap\partial X$ and of the sets $A\cap I$, where $I$ runs over the interval components. Since null sets form a $\sigma$-ideal, the claim follows.
\end{proof}

The same argument applies to any countable cover of $X$ by sets on which $h$ preserves null sets in both directions. We shall often use it in the following form, which we call \emph{finite closed pasting}: if $X=X_1\cup\dots\cup X_k$ and $Y=Y_1\cup\dots\cup Y_k$ with closed sets $X_i,Y_i$, and $h\colon X\to Y$ is a bijection whose restrictions $h|_{X_i}$ are measure-class homeomorphisms onto $Y_i$, then $h$ is a measure-class homeomorphism. Indeed, $h$ is continuous by the pasting lemma, and hence a homeomorphism by compactness.

\begin{theorem}[The null-boundary class]\label{thm:null-class}
Any two Cantorvals with null boundaries are measure-class homeomorphic. The homeomorphism can be chosen as the restriction of an increasing homeomorphism between their convex hulls that is affine on every cell of their boundaries.
\end{theorem}

\begin{proof}
Choose a color-preserving order isomorphism between the cells of the two boundaries, as above, and let $H$ be the homeomorphism given by Lemma~\ref{lem:lift}. Every subset of either boundary is null, so $H$ restricted to the Cantorval is a measure-class homeomorphism by Lemma~\ref{lem:pieces}.
\end{proof}

\subsection{Transport on everywhere fat boundaries}
For ordered Cantor spaces, Akin's measure transport theorem provides order homeomorphisms with arbitrarily small multiplicative distortion between atomless probability measures of full support \cite[Theorem~2.12]{Akin}. The following construction is carried out in the distribution coordinates of Lemma~\ref{lem:distribution}; in addition, it preserves the two colors of cells, which is needed to extend the map to a Cantorval.

\begin{lemma}[Colored interpolation]\label{lem:interpolation}
Let $A_1,A_2,B_1,B_2\subset(0,1)$ be countable dense sets with $A_1\cap A_2=B_1\cap B_2=\varnothing$. For every $L>1$ there is an increasing homeomorphism $\phi\colon[0,1]\to[0,1]$ such that $\phi(A_i)=B_i$ for $i=1,2$ and
\[
L^{-1}|x-y|\le|\phi(x)-\phi(y)|\le L|x-y|\qquad(x,y\in[0,1]).
\]
\end{lemma}

\begin{proof}
Write $A=A_1\cup A_2$ and $B=B_1\cup B_2$, and let $\kappa_A\colon A\to\{1,2\}$ and $\kappa_B\colon B\to\{1,2\}$ be the color maps, so that $a\in A_{\kappa_A(a)}$ and $b\in B_{\kappa_B(b)}$. Choose numbers $0<\epsilon_k<1/4$ with
\[
\prod_k(1-\epsilon_k)>L^{-1},\qquad\prod_k(1+\epsilon_k)<L.
\]
Fix enumerations of $A$ and $B$. We construct finite increasing color-preserving bijections $p_k\colon D_k\to E_k$, where $D_k\subset A\cup\{0,1\}$ and $E_k\subset B\cup\{0,1\}$, each extending the previous one, starting with $p_0(0)=0$ and $p_0(1)=1$. Let $g_k\colon[0,1]\to[0,1]$ be the piecewise affine interpolation of $p_k$. At step $k$ one point is added to $D_{k-1}$, which splits one interval of affinity of $g_{k-1}$ into two. We arrange that each of the two new slopes equals the old slope multiplied by a number in $[1-\epsilon_k,1+\epsilon_k]$.

At an odd step $k$, take the first element $a$ of $A\setminus D_{k-1}$, let $u<a<v$ be its neighbors in $D_{k-1}$, and let $s$ be the slope of $g_{k-1}$ on $[u,v]$. By density we can choose $b\in B_{\kappa_A(a)}\setminus E_{k-1}$ with
\[
|b-g_{k-1}(a)|<\epsilon_k s\min\{a-u,\,v-a\},
\]
and we put $p_k(a)=b$. The two new slopes are $s+(b-g_{k-1}(a))/(a-u)$ and $s-(b-g_{k-1}(a))/(v-a)$, as required.

At an even step $k$, take the first element $b$ of $B\setminus E_{k-1}$ and put $a_*=g_{k-1}^{-1}(b)$. Then $a_*\notin D_{k-1}$; let $u<a_*<v$ be its neighbors in $D_{k-1}$, and let $s$ be the slope of $g_{k-1}$ on $[u,v]$. Choose $a\in A_{\kappa_B(b)}\setminus D_{k-1}$ with
\[
|a-a_*|<\tfrac12\epsilon_k\min\{a_*-u,\,v-a_*\},
\]
and put $p_k(a)=b$. The new slopes are $s(a_*-u)/(a-u)$ and $s(v-a_*)/(v-a)$. The factors lie between $(1+\epsilon_k/2)^{-1}$ and $(1-\epsilon_k/2)^{-1}$, hence in $[1-\epsilon_k,1+\epsilon_k]$.

Every element of $A$ and of $B$ is eventually used. Hence the union of the maps $p_k$ is an increasing color-preserving bijection from $A\cup\{0,1\}$ onto $B\cup\{0,1\}$. All slopes of all the maps $g_k$ lie between $\prod_k(1-\epsilon_k)$ and $\prod_k(1+\epsilon_k)$, so this bijection satisfies the two-sided Lipschitz inequality on the dense set $A\cup\{0,1\}$. Its extension by continuity is the required map $\phi$.
\end{proof}

\begin{lemma}[Boundary distribution coordinates]\label{lem:distribution}
Let $K\subset\R$ be an everywhere fat Cantor set, put $a=\min K$ and $b=\max K$, and define
\[
\mu=\frac{\lambda|_K}{\lambda(K)},\qquad F(x)=\mu([a,x])\quad(a\le x\le b).
\]
Let $P$ be the countable set of endpoints of cells of $K$. Then $F$ is continuous and nondecreasing, $F(K)=[0,1]$, and $F$ is constant on the closure of every cell; distinct cells have distinct values, in the same order as the cells. Moreover, $F_*\mu$ is Lebesgue measure $\ell$ on $[0,1]$, and for every Borel set $C\subset K$,
\[
F(C\setminus P)\ \text{is Borel},\qquad\ell(F(C\setminus P))=\mu(C).
\]
\end{lemma}

\begin{proof}
Since $\mu$ has no atoms, $F$ is continuous, with $F(a)=0$ and $F(b)=1$. A cell has $\mu$-measure zero, so $F$ is constant on its closure; in particular $F(K)=F([a,b])=[0,1]$. If $x<y$ are points of $K$ and $(x,y)$ meets $K$, then $\mu((x,y))>0$, because $\mu$ has full support. Thus two distinct points of $K$ have the same value of $F$ if and only if they are the two endpoints of one cell. This proves the assertion about the order of the values and shows that $F$ is injective on $K\setminus P$.

For $t\in[0,1]$, let $x_t$ be the largest point with $F(x_t)=t$. Then $\mu(F^{-1}([0,t]))=\mu([a,x_t])=t$, which proves $F_*\mu=\ell$. The set $F(C\setminus P)$ is Borel by the Lusin--Souslin theorem \cite[Theorem~15.1]{Kechris}, since $F$ is continuous and injective on the Borel set $C\setminus P$. By the description of the fibers of $F$ given above, $F^{-1}(F(C\setminus P))\cap K=C\setminus P$. Since $\mu$ is concentrated on $K$ and $P$ is countable, we obtain $\ell(F(C\setminus P))=\mu(F^{-1}(F(C\setminus P)))=\mu(C\setminus P)=\mu(C)$.
\end{proof}

\begin{theorem}[The everywhere fat class]\label{thm:fat-class}
Let $X,Y$ be Cantorvals with everywhere fat boundaries $K=\partial X$ and $K'=\partial Y$, and with normalized boundary measures $\mu=\lambda|_K/\lambda(K)$ and $\mu'=\lambda|_{K'}/\lambda(K')$. For every $L>1$ there is an increasing homeomorphism $H$ between their convex hulls, affine on every cell of $K$, such that $H(X)=Y$ and
\[
L^{-1}\mu'(C)\le\mu(H^{-1}(C))\le L\mu'(C)\qquad(C\subset K'\ \text{Borel}).
\]
In particular, $X\mc Y$.
\end{theorem}

\begin{proof}
Let $F$ and $F'$ be the functions given by Lemma~\ref{lem:distribution} for $K$ and $K'$. Assign to each cell $c$ of $K$ the value $d(c)=F(c)\in(0,1)$; this is an order isomorphism from the cells of $K$ onto its image. The images of the cells of interval color and of gap color are both dense in $(0,1)$: the preimage under $F$ of an open subinterval of $(0,1)$ is an open interval meeting $K$, and by Lemma~\ref{lem:cells} it contains cells of both colors. Define $d'$ in the same way for $K'$.

Lemma~\ref{lem:interpolation}, applied to these two pairs of dense sets, gives a map $\phi$, and $\sigma=(d')^{-1}\circ\phi\circ d$ is a color-preserving order isomorphism from the cells of $K$ onto the cells of $K'$. Let $H$ be the homeomorphism given by Lemma~\ref{lem:lift}. Then $F'\circ H=\phi\circ F$ on the closure of every cell, and hence everywhere by continuity. If $P$ and $P'$ are the sets of endpoints of cells of $K$ and $K'$, then $H(P)=P'$, and Lemma~\ref{lem:distribution} gives, for every Borel set $C\subset K'$,
\[
\mu(H^{-1}(C))=\ell\bigl(F(H^{-1}(C)\setminus P)\bigr)=\ell\bigl(\phi^{-1}(F'(C\setminus P'))\bigr).
\]
By the Lipschitz bounds for $\phi^{-1}$, this quantity lies between $L^{-1}\ell(F'(C\setminus P'))=L^{-1}\mu'(C)$ and $L\mu'(C)$. Since every null set is contained in a Borel null set and $H$ maps Borel sets to Borel sets, these bounds show that $H$ maps null subsets of $K$ onto null subsets of $K'$ and vice versa. As $H$ is affine on interval components, Lemma~\ref{lem:pieces} completes the proof.
\end{proof}

\section{Boundary support of achievement sets}\label{sec:support}

Let $x=(x_k)_{k\ge1}$ be a summable sequence of positive numbers, and let
\[
E=E(x)=\Bigl\{\sum_{k=1}^{\infty}\varepsilon_kx_k:\varepsilon_k\in\{0,1\}\Bigr\}.
\]
Rearranging the terms does not change $E$, so throughout this section we assume that $x_1\ge x_2\ge\dots>0$. For $n\ge0$ put
\[
r_n=\sum_{k>n}x_k,\qquad E_n=\Bigl\{\sum_{k>n}\varepsilon_kx_k:\varepsilon_k\in\{0,1\}\Bigr\},
\]
and let $F_n$ be the set of subsums of $x_1,\dots,x_n$, with $F_0=\{0\}$. Thus $E_n$ is the achievement set of the tail $(x_k)_{k>n}$; it is compact, being a continuous image of $\{0,1\}^{\{k>n\}}$. We shall use the relations
\begin{equation}\label{eq:basic}
E=F_n+E_n,\qquad E_n\subset[0,r_n],\qquad E_n=r_n-E_n,\qquad r_n\longrightarrow0.
\end{equation}
All boundaries and supports in this section are taken in $\R$.

\begin{lemma}\label{lem:tails}
Suppose that $\lambda(\partial E)>0$. Then, for every $n\ge0$,
\[
\lambda(\partial E_n)>0\qquad\text{and}\qquad\{0,r_n\}\subset\supp(\lambda|_{\partial E_n}).
\]
\end{lemma}

\begin{proof}
The boundary of a finite union of closed sets is contained in the union of their boundaries. Hence, by~\eqref{eq:basic},
\begin{equation}\label{eq:boundary-cover}
\partial E\subset\bigcup_{f\in F_n}(f+\partial E_n),\qquad\lambda(\partial E)\le|F_n|\,\lambda(\partial E_n),
\end{equation}
which proves the first assertion.

Since $\partial E$ has positive measure, it is infinite, so $E$ is not a finite union of intervals. By Kakeya's theorem~\cite{Kakeya} (see also~\cite{BFP}), the achievement set of a nonincreasing sequence is a finite union of closed intervals if and only if $x_m\le r_m$ for all but finitely many $m$. Hence $x_m>r_m$ for infinitely many indices $m$.

Fix $n$ and choose $m>n$ with $x_m>r_m$. Every nonzero subsum of $x_{n+1},\dots,x_m$ is at least $x_m$, whereas $E_m\subset[0,r_m]$. Since $E_n$ is the union of the translates of $E_m$ by these subsums, we obtain
\begin{equation}\label{eq:tail-agree}
E_n\cap(-\infty,x_m)=E_m.
\end{equation}
In particular, $E_n$ and $E_m$ coincide on the open set $(-\infty,x_m)$, which contains $[0,r_m]$, so $\partial E_m\subset\partial E_n$. For every $\varepsilon>0$ we may choose such an $m$ with $r_m<\varepsilon$, and then, by the first assertion applied to $m$,
\[
\lambda\bigl(\partial E_n\cap(-\varepsilon,\varepsilon)\bigr)\ge\lambda(\partial E_m)>0.
\]
Hence $0\in\supp(\lambda|_{\partial E_n})$. By~\eqref{eq:basic}, the reflection $t\mapsto r_n-t$ maps $E_n$ onto itself; this gives $r_n\in\supp(\lambda|_{\partial E_n})$.
\end{proof}

\begin{lemma}[{\cite[21.8 and 21.11]{BFP}}]\label{lem:hulls}
For $n\ge0$ define
\[
H_n=\bigcup_{f\in F_n}[f,f+r_n].
\]
Then $H_{n+1}\subset H_n$, $E=\bigcap_{n\ge0}H_n$, and $\partial H_n\subset E$. Every gap of $E$ is a gap of $H_n$ for some $n$, that is, a bounded component of $\R\setminus H_n$.
\end{lemma}

These facts are well known; the last assertion is the Second Gap Lemma \cite[21.11]{BFP}.

\begin{theorem}[Support dichotomy]\label{thm:dichotomy}
Let $E$ be the achievement set of a summable sequence of positive numbers. Then
\[
\supp(\lambda|_{\partial E})=
\begin{cases}
\varnothing, & \lambda(\partial E)=0,\\
\partial E, & \lambda(\partial E)>0.
\end{cases}
\]
In particular, the boundary of an achievable Cantorval is either null or everywhere fat.
\end{theorem}

\begin{proof}
Only the case $\lambda(\partial E)>0$ requires proof. Put $S=\supp(\lambda|_{\partial E})$; it is a closed subset of $\partial E$. By Lemma~\ref{lem:tails} with $n=0$, we have $0,r_0\in S$.

Let $(a,b)$ be a gap of $E$. By Lemma~\ref{lem:hulls}, $(a,b)$ is a gap of $H_n$ for some $n$, so that $a=p+r_n$ and $b=q$ for some $p,q\in F_n$. Near $a$ the set $E$ coincides with the translated tail set $p+E_n$, and near $b$ with $q+E_n$: there is $\delta>0$ such that
\begin{align}
E\cap(a-\delta,a+\delta)&=(p+E_n)\cap(a-\delta,a+\delta),\label{eq:NS-left}\\
E\cap(b-\delta,b+\delta)&=(q+E_n)\cap(b-\delta,b+\delta).\label{eq:NS-right}
\end{align}
This is the Nymann--S{\'a}enz lemma; see \cite[21.19]{BFP}. (It also follows directly from the fact that all the intervals $[f,f+r_n]$, $f\in F_n$, have the same length: $[p,p+r_n]$ is the only one with right endpoint $a$, $[q,q+r_n]$ is the only one with left endpoint $b$, and all the others lie at positive distance from $a$ and from $b$, respectively.)

Two sets that coincide on an open set have the same boundary there. Therefore \eqref{eq:NS-left}, \eqref{eq:NS-right} and Lemma~\ref{lem:tails} show that every neighborhood of $a=p+r_n$ and every neighborhood of $b=q+0$ meets $\partial E$ in a set of positive measure. Hence $a,b\in S$.

Finally, the endpoints of the gaps of $E$, together with $0$ and $r_0$, are dense in $\partial E$. Indeed, let $y\in\partial E$ and $\eta>0$. There is a point $z\notin E$ with $|z-y|<\eta$. The component of $\R\setminus E$ containing $z$ is either a gap of $E$ or one of the half-lines $(-\infty,0)$ and $(r_0,\infty)$, and its endpoint nearer to $y$ lies between $y$ and $z$ or coincides with $y$. Thus $S$ contains a dense subset of $\partial E$, and since $S$ is closed, $S=\partial E$.
\end{proof}

\section{Isolated digits and positive boundary measure}\label{sec:positive}

\subsection{Sets of generalized expansions}
Let $D_j\subset[0,\infty)$ be finite sets containing $0$, and let $w_j>0$ satisfy
\[
\sum_{j=1}^{\infty}w_j\max D_j<\infty.
\]
We consider the set
\begin{equation}\label{eq:expansions}
E=\Bigl\{\sum_{j=1}^{\infty}w_jd_j:d_j\in D_j\Bigr\},\qquad R_n=\sum_{j>n}w_j\max D_j.
\end{equation}
Sets of this form are known as sets of generalized expansions; see~\cite{AC}. The set $E$ is compact and $R_n\to0$. If each $D_j$ is the set of subsums of a finite block of positive numbers, then concatenating these blocks, multiplied by the weights $w_j$, shows that $E$ is an achievement set.

For a word $v=(d_1,\dots,d_n)\in D_1\times\dots\times D_n$ put $s_v=\sum_{j=1}^nw_jd_j$. We call $[s_v,s_v+R_n]$ the \emph{hull} of $v$; it contains every point $\sum_jw_jd_j$ of $E$ whose first $n$ digits are $d_1,\dots,d_n$. In particular, $E$ is contained in the union of the hulls of all words of length $n$. Since $w_{n+1}d+R_{n+1}\le R_n$ for every $d\in D_{n+1}$, the hull of every one-letter extension of $v$ is contained in the hull of $v$.

\begin{definition}\label{def:isolated}
A family of nonempty sets $A_j\subset D_j$ is \emph{isolated} if
\begin{equation}\label{eq:isolated}
w_j|a-d|>R_j\qquad(a\in A_j,\ d\in D_j\setminus\{a\}).
\end{equation}
Its associated compact set is
\[
K=\Bigl\{\sum_{j=1}^{\infty}w_ja_j:a_j\in A_j\Bigr\}.
\]
\end{definition}

Condition~\eqref{eq:isolated} says the following. Let two words of length $j$ coincide in their first $j-1$ letters, and let the last letter of one of them be a selected digit $a\in A_j$. Then the hulls of the two words are disjoint, and their distance is at least $w_j|a-d|-R_j>0$, where $d$ is the last letter of the other word.

\begin{theorem}[Isolated digit criterion]\label{thm:isolated}
Let $(A_j)$ be an isolated family. Then every point of $K$ is a one-point connected component of $E$. Moreover,
\begin{equation}\label{eq:isolated-measure}
K\subset\partial E,\qquad\lambda(K)=\lim_{n\to\infty}R_n\prod_{j=1}^n|A_j|.
\end{equation}
\end{theorem}

\begin{proof}
Call a word $u=(a_1,\dots,a_n)$ \emph{selected} if $a_j\in A_j$ for all $j\le n$. Let $I_u=[s_u,s_u+R_n]$ be its hull and put
\[
K_n=\bigcup_{u\in A_1\times\dots\times A_n}I_u.
\]
Since hulls are nested, $K_{n+1}\subset K_n$.

Let $u$ be a selected word and $v\ne u$ any word of length $n$, and let $j$ be the first position at which they differ. By~\eqref{eq:isolated}, the hulls of $(u_1,\dots,u_j)$ and $(v_1,\dots,v_j)$ are at distance at least $w_j|u_j-v_j|-R_j>0$, and the hulls of $u$ and $v$ are contained in them. Hence $I_u$ lies at positive distance from the hull of every other word of length $n$. In particular, the intervals $I_u$ are pairwise disjoint, and $\lambda(K_n)=R_n\prod_{j=1}^n|A_j|$.

Clearly $K\subset K_n$ for every $n$. Conversely, let $x\in\bigcap_nK_n$. For each $n$ there is exactly one selected word $u^{(n)}$ of length $n$ with $x\in I_{u^{(n)}}$, and by nesting and disjointness $u^{(n+1)}$ extends $u^{(n)}$. Since $R_n\to0$, the value of the resulting infinite selected word is $x$. Thus $K=\bigcap_nK_n$, and continuity of measure from above gives the formula in~\eqref{eq:isolated-measure}.

Finally, let $x\in K$ and $n\ge1$, and let $u$ be the selected word of length $n$ with $x\in I_u$. The set $E$ is contained in the union of the hulls of all words of length $n$, and $I_u$ lies at positive distance from all of them except itself. Hence $E\cap I_u$ is open and closed in $E$, and the connected component of $x$ in $E$ is contained in $I_u$; its diameter is at most $R_n$. Letting $n\to\infty$, we see that this component is $\{x\}$. In particular $x\notin\intr E$, so $x\in\partial E$.
\end{proof}

\subsection{Hybrid binary blocks}
Let $a_j\ge2$ be integers, and set
\[
b_j=2^{a_j},\qquad Q_0=1,\qquad Q_j=\prod_{i=1}^jb_i,\qquad t_j=\sum_{i=1}^ja_i.
\]
Note that $Q_j=2^{t_j}$. Consider the sequence obtained by concatenating, for $j=1,2,\dots$, the blocks
\begin{equation}\label{eq:blocks}
\frac{1}{Q_j}\Bigl(b_j-1,\ \frac{b_j}{2},\ \frac{b_j}{4},\ \dots,\ 2\Bigr).
\end{equation}
The $j$-th block has $a_j$ terms, so the first $j$ blocks consist of $t_j$ terms. These are the hybrid series introduced by Pratsiovytyi, Vynnyshyn and Vynnyshyn~\cite{PVV}. Let $E$ denote the achievement set, let $E^{(n)}=E_{t_n}$ be the achievement set of the terms that follow the first $n$ blocks, and let $R_n$ be the sum of these terms. Since the subsums of $2,4,\dots,b_j/2$ are exactly the even integers $0,2,\dots,b_j-2$, the subsums of the $j$-th block, multiplied by $Q_j$, form the set
\begin{equation}\label{eq:digits}
D_j=\{0,2,\dots,b_j-2\}\cup\{b_j-1,b_j+1,\dots,2b_j-3\}.
\end{equation}
Thus $E$ is the set~\eqref{eq:expansions} with these digit sets and the weights $w_j=1/Q_j$, and $R_n$ agrees with the remainder in~\eqref{eq:expansions}.

\begin{lemma}\label{lem:tail-sums}
The sequence~\eqref{eq:blocks} is positive, strictly decreasing and summable, and
\begin{equation}\label{eq:tail-sums}
R_n=\frac{2}{Q_n}-\sum_{j>n}\frac{1}{Q_j},\qquad\frac{5}{3Q_n}\le R_n<\frac{2}{Q_n}.
\end{equation}
\end{lemma}

\begin{proof}
Within a block the terms decrease, since $b_j-1>b_j/2$. At the junction of the $j$-th and the $(j+1)$-st block,
\[
\frac{b_{j+1}-1}{Q_{j+1}}<\frac{b_{j+1}}{Q_{j+1}}=\frac{1}{Q_j}<\frac{2}{Q_j}.
\]
The sum of the $j$-th block is
\[
\frac{2b_j-3}{Q_j}=\frac{2}{Q_{j-1}}-\frac{2}{Q_j}-\frac{1}{Q_j},
\]
and telescoping gives the formula for $R_n$ in~\eqref{eq:tail-sums}. Since $b_j\ge4$, we have $\sum_{j>n}Q_n/Q_j\le\sum_{k\ge1}4^{-k}=1/3$, which gives the bounds in~\eqref{eq:tail-sums}. In particular, $R_0<2$.
\end{proof}

\begin{lemma}\label{lem:hybrid-cantorval}
The achievement set $E$ of every sequence~\eqref{eq:blocks} is a Cantorval \textup{\cite[Theorem~4.2]{PVV}}. Moreover, $E$ contains the interval $[R_0-1,1]$.
\end{lemma}

\begin{proof}
Only the second assertion needs proof. It slightly improves the inclusion $[1-1/b_1,1]\subset E$ of \cite[Theorem~3.3]{PVV}, since $R_0-1=1-\sum_{j\ge1}Q_j^{-1}<1-1/b_1$.

The set $D_j$ contains exactly one representative of each residue class modulo $b_j$. Put
\[
P_n=\Bigl\{\sum_{j=1}^n\frac{d_j}{Q_j}:d_j\in D_j\Bigr\}.
\]
We claim that $Q_nP_n$ contains exactly one representative of each residue class modulo $Q_n$. Since there are $Q_n$ words $(d_1,\dots,d_n)$, it suffices to show that distinct words give distinct residues. Write $Q_np=b_n(Q_{n-1}p')+d_n$, where $p'\in P_{n-1}$ comes from the first $n-1$ digits. The residue of $Q_np$ modulo $b_n$ determines $d_n$, and then the residue of $Q_{n-1}p'$ modulo $Q_{n-1}$ is determined; the claim follows by induction.

Let $\pi\colon\R\to\R/\mathbb{Z}$ be the quotient map. Since $P_n\subset E$, the compact set $\pi(E)$ contains $\pi(Q_n^{-1}\mathbb{Z})$ for every $n$. These grids are dense, so $\pi(E)=\R/\mathbb{Z}$. If $R_0-1<t<1$, then $t$ is the only point of $[0,R_0]\supset E$ that is congruent to $t$ modulo $1$, because $R_0<2$. Hence $t\in E$, and since $E$ is closed, $[R_0-1,1]\subset E$.
\end{proof}

\subsection{The boundary measure}
The set~\eqref{eq:digits} contains exactly one pair of consecutive integers, namely $b_j-2$ and $b_j-1$, and any two other elements of $D_j$ differ by at least $2$. Define
\begin{equation}\label{eq:selected}
A_j=D_j\setminus\{b_j-2,b_j-1\}=\{0,2,\dots,b_j-4\}\cup\{b_j+1,b_j+3,\dots,2b_j-3\}.
\end{equation}
Every $a\in A_j$ differs from every other element $d$ of $D_j$ by at least $2$, so by~\eqref{eq:tail-sums},
\[
\frac{|a-d|}{Q_j}\ge\frac{2}{Q_j}>R_j.
\]
Hence the family $(A_j)$ is isolated in the sense of Definition~\ref{def:isolated}.

\begin{theorem}\label{thm:hybrid}
Let $E$ be the achievement set of the sequence~\eqref{eq:blocks}, and let $K$ be the compact set associated with the family~\eqref{eq:selected}. Then $K$ is a Cantor set contained in $\partial E$, and
\begin{equation}\label{eq:product}
\lambda(K)=2\prod_{j=1}^{\infty}\Bigl(1-\frac{2}{b_j}\Bigr).
\end{equation}
In particular,
\[
\sum_{j=1}^{\infty}\frac{1}{b_j}<\infty\implies\lambda(\partial E)\ge2\prod_{j=1}^{\infty}\Bigl(1-\frac{2}{b_j}\Bigr)>0.
\]
For $a_j=j+1$, the resulting Cantorval $E_*$ satisfies
\begin{equation}\label{eq:E-star}
\lambda(\partial E_*)\ge2\prod_{j=1}^{\infty}(1-2^{-j})=0.5775761901732048\ldots.
\end{equation}
\end{theorem}

\begin{proof}
Each $A_j$ has $b_j-2\ge2$ elements, so every point of $K$ is a limit of other points of $K$, and $K$ is perfect. By Theorem~\ref{thm:isolated}, the points of $K$ are one-point components of $E$; hence $K$ is totally disconnected, and $K\subset\partial E$. Thus $K$ is a Cantor set contained in $\partial E$. Since $|A_j|=b_j-2=b_j(1-2/b_j)$ and $Q_n=b_1\cdots b_n$, formula~\eqref{eq:isolated-measure} becomes
\begin{equation}\label{eq:limit}
\lambda(K)=\lim_{n\to\infty}(Q_nR_n)P_n,\qquad P_n=\prod_{j=1}^n\Bigl(1-\frac{2}{b_j}\Bigr).
\end{equation}
The partial products $P_n$ decrease to $P=\prod_{j\ge1}(1-2/b_j)\ge0$, and by~\eqref{eq:tail-sums} we have $5/3\le Q_nR_n<2$ for all $n$. We distinguish two cases.

If $P=0$, then $0\le Q_nR_nP_n<2P_n\to0$, so $\lambda(K)=0=2P$, which is~\eqref{eq:product}.

If $P>0$, then $\sum_j2/b_j\le-\sum_j\log(1-2/b_j)=-\log P<\infty$, because $-\log(1-u)\ge u$ for $0\le u<1$; in particular, $b_j\to\infty$. By~\eqref{eq:tail-sums},
\[
0<2-Q_nR_n=\sum_{j>n}\frac{Q_n}{Q_j}=\frac{1}{b_{n+1}}\Bigl(1+\frac{1}{b_{n+2}}+\frac{1}{b_{n+2}b_{n+3}}+\cdots\Bigr)\le\frac{4}{3b_{n+1}}\longrightarrow0.
\]
Thus $Q_nR_n\to2$, and~\eqref{eq:limit} gives $\lambda(K)=2P$, which is again~\eqref{eq:product}.

Conversely, if $\sum_j1/b_j<\infty$, then $P>0$, since $\log(1-u)\ge-2u$ for $0\le u\le1/2$ and $2/b_j\le1/2$. Together with $K\subset\partial E$, this proves the second assertion. Finally, for $a_j=j+1$ we have $b_j=2^{j+1}$ and $2/b_j=2^{-j}$, and~\eqref{eq:E-star} follows.
\end{proof}

For $a_j=j+1$ one has $Q_j=2^{j(j+3)/2}$, and the sequence begins
\[
\frac34,\ \frac12;\quad\frac{7}{32},\ \frac18,\ \frac{1}{16};\quad\frac{15}{512},\ \frac{1}{64},\ \frac{1}{128},\ \frac{1}{256};\quad\dots.
\]
This is the series of \cite[Example~2.5]{PVV}. Formula~\eqref{eq:product} gives the exact measure of the Cantor set $K\subset\partial E_*$ constructed above, and the sum of the series is $R_0=1.7167348786896923\ldots$.

\begin{corollary}\label{cor:close-to-one}
For every $\varepsilon>0$ there is an achievable Cantorval $\widetilde{E}\subset[0,1]$ with $\min\widetilde{E}=0$ and $\max\widetilde{E}=1$ such that
\[
\lambda(\partial\widetilde{E})>1-\varepsilon.
\]
\end{corollary}

\begin{proof}
Take $a_j=j+M$ with an integer $M\ge2$. Then
\[
\prod_{j\ge1}\Bigl(1-\frac{2}{b_j}\Bigr)\ge1-\sum_{j\ge1}\frac{2}{b_j}=1-2^{1-M}.
\]
Let $\widetilde{E}=E/R_0$, the achievement set of the sequence~\eqref{eq:blocks} divided by its sum. Theorem~\ref{thm:hybrid} and $R_0<2$ give $\lambda(\partial\widetilde{E})=\lambda(\partial E)/R_0>1-2^{1-M}$. It remains to choose $M$ with $2^{1-M}<\varepsilon$.
\end{proof}

\section{The two achievable measure classes}\label{sec:two-classes}

Let
\[
\Gamma=E\Bigl(\frac{3}{4^j},\frac{2}{4^j}:j\ge1\Bigr)
\]
be the Guthrie--Nymann Cantorval, that is, the achievement set of the sequence $3/4$, $2/4$, $3/16$, $2/16$, \dots, and let $E_*$ be the hybrid Cantorval of Theorem~\ref{thm:hybrid}. The boundary of $\Gamma$ is Lebesgue null~\cite{NPW}.

\begin{theorem}\label{thm:two-classes}
Every achievable Cantorval $X$ belongs to exactly one of the following two measure-class homeomorphism classes:
\begin{align*}
\mathsf{N}&:\quad X\mc\Gamma,\quad\lambda(\partial X)=0,\\
\mathsf{F}&:\quad X\mc E_*,\quad\lambda(\partial X)>0.
\end{align*}
If $X$ and $Y$ belong to the same class, a measure-class homeomorphism between them can be chosen as the restriction of an increasing homeomorphism between their convex hulls that is affine on every interval component and on every gap. If $X,Y\in\mathsf{F}$, then for every $L>1$ such a map $H$ can be chosen to satisfy
\begin{equation}\label{eq:distortion}
L^{-1}\widehat{\mu}_Y(B)\le\widehat{\mu}_X(H^{-1}(B))\le L\widehat{\mu}_Y(B)\qquad(B\subset\partial Y\ \text{Borel}),
\end{equation}
where $\widehat{\mu}_X=\lambda|_{\partial X}/\lambda(\partial X)$.
\end{theorem}

\begin{proof}
By Theorem~\ref{thm:dichotomy}, either $S(X)=\varnothing$ or $S(X)=\partial X$. In the first case, $X\mc\Gamma$ by Theorem~\ref{thm:null-class}. In the second case the boundary is everywhere fat, and Theorem~\ref{thm:fat-class} gives $X\mc E_*$, as well as~\eqref{eq:distortion} for any two members of $\mathsf{F}$. Measure-class homeomorphisms map boundaries onto boundaries (Lemma~\ref{lem:homeo-invariance}) and preserve null sets, so the two classes are disjoint. Both are nonempty, as $\Gamma$ and $E_*$ show.
\end{proof}

\subsection{Tail interior density and standardness}\label{ss:standard}
For a nonincreasing summable sequence of positive numbers we keep the notation $E_n$, $r_n$, $F_n$ of Section~\ref{sec:support}, and put
\[
\theta_n=\frac{\lambda(\intr E_n)}{r_n},\qquad\ell_n=\sup\{|J|:J\ \text{is an interval component of}\ E_n\},
\]
with $\ell_n=0$ if $E_n$ has no interval components.

\begin{proposition}\label{prop:interior-density}
If $\limsup_n\theta_n>0$, then $\lambda(\partial E)=0$. Consequently, if $\lambda(\partial E)>0$, then
\begin{equation}\label{eq:standard}
\lim_{n\to\infty}\frac{\lambda(\intr E_n)}{r_n}=\lim_{n\to\infty}\frac{\ell_n}{r_n}=0.
\end{equation}
\end{proposition}

\begin{proof}
Choose $\delta>0$ and infinitely many $n$ such that $\lambda(\intr E_n)\ge2\delta r_n$. For $y\in\partial E$ and such an $n$, choose $f\in F_n$ with $y\in f+E_n$; this is possible by~\eqref{eq:basic}. Then
\[
f+E_n\subset[y-r_n,y+r_n],\qquad f+\intr E_n\subset\intr E,
\]
and therefore
\[
\frac{\lambda(\partial E\cap[y-r_n,y+r_n])}{2r_n}\le1-\delta.
\]
Since $r_n\to0$, no point of $\partial E$ is a density point of $\partial E$, and the Lebesgue density theorem gives $\lambda(\partial E)=0$. The last assertion follows by contraposition, since $0\le\ell_n\le\lambda(\intr E_n)$.
\end{proof}

Following Nowakowski and Prus-Wi{\'s}niowski~\cite{NPW}, the achievable Cantorval $E(x)$ generated by a nonincreasing sequence $x$ is called \emph{standard} if
\[
\limsup_{n\to\infty}\frac{\ell_n}{r_n}>0.
\]
Proposition~\ref{prop:interior-density} contains the theorem that standard Cantorvals have null boundaries \cite[Theorem~5]{NPW}.

\begin{corollary}\label{cor:counterexample}
The Cantorval $E_*$ has a boundary of positive Lebesgue measure, and none of its nonincreasing generating sequences is standard. This answers Problem~21 of the preprint version of~\cite{GPW} in the affirmative, and gives a counterexample to \cite[Conjecture~11]{NPW}, which also appears as Problem~22 in the same survey.
\end{corollary}

\begin{proof}
Apply Theorem~\ref{thm:hybrid}, and Proposition~\ref{prop:interior-density} to any nonincreasing sequence of positive numbers whose achievement set is $E_*$.
\end{proof}

\begin{proposition}\label{prop:hybrid-standard}
For the hybrid sequence~\eqref{eq:blocks},
\[
b_j\longrightarrow\infty\implies\lim_{n\to\infty}\frac{\ell_n}{r_n}=0.
\]
\end{proposition}

\begin{proof}
Let $L_j$ be the maximal length of a connected component of $E^{(j)}$. By~\eqref{eq:digits},
\[
E^{(j-1)}=\bigcup_{d\in D_j}\Bigl(\frac{d}{Q_j}+E^{(j)}\Bigr).
\]
Each of these translates lies in its hull $[d/Q_j,d/Q_j+R_j]$, of length $R_j<2/Q_j$. Consecutive elements of $D_j$ differ by $2$, except for the pair $b_j-2,b_j-1$. Hence the hulls are pairwise disjoint, except for the two hulls corresponding to $b_j-2$ and $b_j-1$, and every connected component of $E^{(j-1)}$ lies in a single hull or in the union of these two hulls. It follows that
\[
L_{j-1}\le R_j+\frac{1}{Q_j}<\frac{3}{Q_j},\qquad\frac{L_{j-1}}{R_{j-1}}\le\frac{3}{2b_j-3},
\]
where we used that $R_{j-1}$ is at least the sum $(2b_j-3)/Q_j$ of the $j$-th block.

Now let $n$ be an index in the $j$-th block, that is, $t_{j-1}<n\le t_j$. The terms of the $j$-th block that come after the $n$-th term are of the form $2^i/Q_j$, so their distinct subsums are even multiples of $1/Q_j$ and differ by at least $2/Q_j>R_j$. Hence $E_n$ is a union of translates of $E^{(j)}$ whose hulls lie at positive distance from each other. Every component of $E_n$ is therefore a translate of a component of $E^{(j)}$, so $\ell_n=L_j$, while $r_n\ge R_j$. This gives
\[
\frac{\ell_n}{r_n}\le\frac{L_j}{R_j}\le\frac{3}{2b_{j+1}-3}.
\]
Every index $n\ge1$ lies in some block, and the right-hand side tends to $0$ as $j\to\infty$.
\end{proof}

\section{Clopen supports and interface invariants}\label{sec:clopen}

The measure-class framework also applies to Cantorvals whose boundary support is a nonempty proper subset of the boundary. Throughout this section $X$ is a Cantorval, and
\[
K=\partial X,\qquad S=S(X),\qquad R=R(X)=K\setminus S.
\]
The set $R$ is null, by the remark after Definition~\ref{def:support}. The set $S$ is closed and, if nonempty, perfect: an isolated point of the support of an atomless measure would have a neighborhood of measure zero.

\subsection{Finite decompositions}

\begin{lemma}\label{lem:clopen-pieces}
Every nonempty clopen subset $C$ of a Cantorval $X$ is a Cantorval, and $\partial C=C\cap\partial X$. If the interiors of finitely many interval components are removed from $X$, the remaining set is a Cantorval with boundary $\partial X$.
\end{lemma}

\begin{proof}
If $C\ne X$, the compact sets $C$ and $X\setminus C$ are at positive distance from each other, so $C$ coincides with $X$ in a neighborhood of each of its points. Hence $C$ is regular closed, $\partial C=C\cap\partial X$, and this set is perfect and nowhere dense, that is, a Cantor set.

Let $Y$ be obtained by removing the interiors of finitely many interval components. It is compact, and its interior is the union of the interiors of the remaining interval components. Let $p$ be an endpoint of a removed component, say its left endpoint. Since $\partial X$ is perfect and the component lies to the right of $p$, points of $\partial X$ accumulate at $p$ from the left; since $X$ is regular closed, so do interior points of $X$, and these can be taken in components different from the finitely many removed ones. The same holds near every other point of $\partial X$. Hence $Y=\cl(\intr Y)$ and $\partial Y=\partial X$.
\end{proof}

A point of a Cantorval $X$ is \emph{exposed} if it is an endpoint of a component of $\R\setminus X$, that is, an endpoint of a gap of $X$ or one of the extreme points $\min X$ and $\max X$.

\begin{lemma}[Finite exposed-point matching]\label{lem:matching}
Let $Y,Y'$ be Cantorvals whose boundaries are both null or both everywhere fat. Given distinct exposed points $p_1,\dots,p_t$ of $Y$ and distinct exposed points $p'_1,\dots,p'_t$ of $Y'$, there is a measure-class homeomorphism $g\colon Y\to Y'$, affine on every interval component, such that $g(p_i)=p'_i$ for all $i$.
\end{lemma}

\begin{proof}
Cutting $Y$ at finitely many gaps, we partition it into finitely many clopen pieces, which are Cantorvals by Lemma~\ref{lem:clopen-pieces}, in such a way that every marked point $p_i$ is an extreme point of its piece and no piece contains two marked points. For this, cut at the gap adjacent to each marked point that is not an extreme point of $Y$, and at further gaps separating the marked points. Cutting a marked piece once more, at a gap away from its marked point, produces an unmarked piece. Proceed in the same way with $Y'$, arrange that both partitions contain an unmarked piece, and subdivide unmarked pieces until the two partitions have the same number of pieces. The boundary of each piece is a relatively open subset of the original boundary, so it is null or everywhere fat together with it.

Now match the marked pieces according to the prescribed pairs, and the unmarked pieces arbitrarily. For a pair of marked pieces, if $p_i$ and $p'_i$ are both left or both right extreme points of their pieces, apply Theorem~\ref{thm:null-class} or Theorem~\ref{thm:fat-class} to the two pieces; otherwise apply it to the first piece and the mirror image of the second, and then reflect back. This gives a measure-class homeomorphism between the two pieces that sends $p_i$ to $p'_i$. Between unmarked pieces use an increasing map given by the same theorem. These finitely many maps paste to the required homeomorphism by finite closed pasting.
\end{proof}

Assume now that $S$ is clopen in $K$. An interval component of $X$ with one endpoint in $S$ and the other in $R$ is called a \emph{bridge}. We write $m(X)$ for the number of bridges.

\begin{remark*}
In the theory of Cantor sets the word \emph{bridge} classically has a different meaning, which goes back to Newhouse's notion of thickness~\cite{Newhouse}. For a gap $U$ of a Cantor set and an endpoint $u$ of $U$, the bridge at $u$ is the maximal interval with endpoint $u$ that contains no point of a gap of length at least $|U|$, the unbounded components of the complement included. Thus its length is the distance from $u$ to the nearest such gap on the side opposite to $U$. Our bridges are unrelated to this notion; the name refers to the fact that they join the two parts $S$ and $R$ of the boundary.
\end{remark*}

\begin{lemma}[Bridge decomposition]\label{lem:bridges}
Suppose that $S$ and $R$ are nonempty and clopen in $K$. Then there are only finitely many bridges $B_1,\dots,B_m$. The sets
\begin{align*}
X_S&=S\cup\bigcup\{I: I\ \text{is an interval component of}\ X,\ \partial I\subset S\},\\
X_R&=R\cup\bigcup\{I: I\ \text{is an interval component of}\ X,\ \partial I\subset R\}
\end{align*}
are disjoint Cantorvals with everywhere fat and null boundaries, respectively, and
\[
\partial X_S=S,\qquad\partial X_R=R,\qquad X=X_S\cup X_R\cup B_1\cup\dots\cup B_m.
\]
Each bridge meets each of $X_S$ and $X_R$ in exactly one point, and this point is exposed in the respective set.
\end{lemma}

\begin{proof}
Let $\delta=\dist(S,R)>0$. Every bridge has length at least $\delta$, so there are only finitely many of them. To see that $X_S$ is closed, consider a convergent sequence of points of $X_S$. If infinitely many of them lie in $S$ or in finitely many interval components, the limit clearly belongs to $X_S$. Otherwise, after passing to a subsequence, the points lie in distinct interval components with endpoints in $S$. The lengths of these components tend to $0$, so their endpoints converge to the same limit, which belongs to the closed set $S$. The same argument applies to $X_R$.

Since $S$ and $R$ are clopen subsets of the Cantor set $K$, they are perfect. In a sufficiently small neighborhood of a point of $S$ all points of $K$ belong to $S$, and Lemma~\ref{lem:cells} provides there cells of interval color, that is, interval components of $X$ with both endpoints in $S$. Hence $X_S$ is regular closed with boundary $S$. Similarly, $X_R$ is a Cantorval with boundary $R$. Since $K\setminus S$ is null, $\lambda|_S=\lambda|_K$, and its support is $S$; thus $\partial X_S$ is everywhere fat, while $\partial X_R=R$ is null. The remaining assertions follow by sorting the interval components of $X$ according to the positions of their endpoints; after a bridge is removed, its endpoint becomes an exposed point of the corresponding part.
\end{proof}

\begin{theorem}[Classification with clopen support]\label{thm:clopen}
The Cantorvals $X$ for which $S(X)$ is clopen in $\partial X$ fall into precisely the following measure-class homeomorphism types:
\[
\mathsf{N},\qquad\mathsf{F},\qquad\mathsf{M}_m\quad(m=0,1,2,\dots).
\]
Here $\mathsf{N}$ consists of the Cantorvals with null boundary, $\mathsf{F}$ of those with everywhere fat boundary, and $\mathsf{M}_m$ of those for which $S(X)$ and $R(X)$ are both nonempty and there are exactly $m$ bridges. Every listed type occurs.
\end{theorem}

\begin{proof}
By Lemmas~\ref{lem:supp} and~\ref{lem:homeo-invariance}, a measure-class homeomorphism preserves $S$, $R$, the interval components and their endpoints, and hence the number of bridges. The pure types $\mathsf{N}$ and $\mathsf{F}$ are single classes by Theorems~\ref{thm:null-class} and~\ref{thm:fat-class}.

Consider two Cantorvals of type $\mathsf{M}_m$, and choose a bijection between their bridges. In the decompositions of Lemma~\ref{lem:bridges}, apply Lemma~\ref{lem:matching} to the two fat parts and to the two null parts, prescribing the correspondence between the endpoints of bridges induced by this bijection. Extend affinely over every bridge. These maps agree where the finitely many closed pieces meet, and they paste to a bijection, which is a measure-class homeomorphism by finite closed pasting.

To see that every type occurs, choose Cantorvals $\Phi$ and $\Psi$ with convex hull $[0,1]$ of types $\mathsf{F}$ and $\mathsf{N}$, respectively; for instance, rescaled copies of $E_*$ and $\Gamma$ (Theorems~\ref{thm:dichotomy} and~\ref{thm:hybrid}, and~\cite{NPW}). They realize the types $\mathsf{F}$ and $\mathsf{N}$. For $m\ge1$, place copies of $\Phi$ and $\Psi$ alternately on the intervals $[3i,3i+1]$, $0\le i\le m$, and adjoin the intervals $[3i+1,3i+3]$, $0\le i<m$. The extreme points of each copy are accumulation points of its gaps from the inner side, so they remain boundary points after the intervals are adjoined, and no adjoined interval merges with an interval component of a copy. The resulting compact set is regular closed, its boundary is the union of the boundaries of the copies, hence a Cantor set, and it has exactly $m$ bridges. Its boundary support is the union of the boundaries of the copies of $\Phi$, which is clopen in the boundary. For $m=0$, take two disjoint copies, one of $\Phi$ and one of $\Psi$.
\end{proof}

By Theorem~\ref{thm:dichotomy}, no Cantorval of type $\mathsf{M}_m$ is an achievement set.

\subsection{The interface when the support is not open}
Define the \emph{support interface} of $X$ by
\[
\mathcal{I}(X)=S(X)\cap\cl_{\partial X}\bigl(\partial X\setminus S(X)\bigr).
\]
By Lemmas~\ref{lem:supp} and~\ref{lem:homeo-invariance}, a measure-class homeomorphism maps the support interface of one Cantorval onto that of the other. Consequently, the nullity of $\mathcal{I}(X)$ is a measure-class invariant.

\begin{proposition}\label{prop:interface}
There exist Cantorvals $X_0,X_1$ and a homeomorphism $h\colon X_0\to X_1$ with $h(S(X_0))=S(X_1)$ such that $X_0$ and $X_1$ are not measure-class homeomorphic.
\end{proposition}

\begin{proof}
Fix Cantorvals $\Phi$ and $\Psi$ with convex hull $[0,1]$ whose boundaries are everywhere fat and null, respectively. Let $T\subset\R$ be a Cantor set, and partition its cells into two families, each of which is dense in the order of cells. In the closed middle third of each cell $J$ place an affine copy $P_J$ of $\Phi$ or of $\Psi$, according to the family of $J$, with that middle third as its convex hull. Put
\[
X(T)=T\cup\bigcup_JP_J.
\]
For every $\varepsilon>0$ only finitely many cells have length greater than $\varepsilon$, so every limit of points taken from infinitely many distinct pieces lies in $T$; hence $X(T)$ is compact. Every point of $T$ is approached by cells of both families, and hence by interior points of the inserted pieces. It follows that $X(T)$ is regular closed and
\begin{gather*}
\partial X(T)=T\cup\bigcup_J\partial P_J,\qquad S(X(T))=T\cup\bigcup_{J\ \mathrm{fat}}\partial P_J,\\
\partial X(T)\setminus S(X(T))=\bigcup_{J\ \mathrm{null}}\partial P_J.
\end{gather*}
For the formula for the support, note that every neighborhood of a point of $T$ contains a piece of the fat family, while each piece of the null family has a neighborhood disjoint from $T$ and from all other pieces. These formulas also show that the boundary is perfect and nowhere dense, so that $X(T)$ is a Cantorval, and that $\mathcal{I}(X(T))=T$.

Now let $T_0$ be a null Cantor set and $T_1$ an everywhere fat Cantor set. Choose an order isomorphism between their cells, and let $h$ be the increasing homeomorphism between their convex hulls given by Lemma~\ref{lem:lift}; it is affine on every cell and maps $T_0$ onto $T_1$. Construct $X_0=X(T_0)$ and put $X_1=h(X_0)$. Since $h$ is affine on each cell, $X_1=X(T_1)$ for the transported partition of the cells of $T_1$, with inserted pieces of the same types. The formulas above give $h(S(X_0))=S(X_1)$ and $\mathcal{I}(X_i)=T_i$ for $i=0,1$. A measure-class homeomorphism from $X_0$ onto $X_1$ would map the null set $T_0=\mathcal{I}(X_0)$ onto the set $T_1=\mathcal{I}(X_1)$ of positive measure, which is impossible.
\end{proof}

\section{The group of measure-class homeomorphisms}\label{sec:groups}

For a Cantorval $X$, let $\Homeo(X)$ carry the metric
\[
\rho(h,g)=\sup_{x\in X}|h(x)-g(x)|+\sup_{x\in X}|h^{-1}(x)-g^{-1}(x)|.
\]
With this metric $\Homeo(X)$ is a Polish group~\cite{Kechris}. We write $\MC(X)$ for its subgroup of measure-class homeomorphisms, and $\mathcal{C}(X)$ for the countable family of interval components of $X$.

\subsection{Generic singularity on interval components}

\begin{theorem}\label{thm:meager}
For every Cantorval $X$, the subgroup $\MC(X)$ is meager in $\Homeo(X)$. More precisely, for a comeager set of $h\in\Homeo(X)$ and every $I\in\mathcal{C}(X)$, there is a Borel null set $N\subset I$ such that
\[
\lambda(h(N))\ge\tfrac12\lambda(h(I)).
\]
\end{theorem}

\begin{proof}
Fix $I\in\mathcal{C}(X)$. For $J\in\mathcal{C}(X)$ put
\[
\mathcal{O}_J=\{h\in\Homeo(X):h(I)=J\}.
\]
By Lemma~\ref{lem:homeo-invariance} these sets partition $\Homeo(X)$. They are open: if $c$ is the midpoint of $I$ and $h\in\mathcal{O}_J$, then every $g$ sufficiently close to $h$ satisfies $g(c)\in\intr J$, which forces $g(I)=J$. Hence each $\mathcal{O}_J$ is also closed, its complement being the union of the other members of the partition.

For $k\ge1$, let $\mathcal{B}_{J,k}$ consist of the maps $h\in\mathcal{O}_J$ such that
\[
\lambda(h(U))\le|J|/2
\]
for every finite union $U$ of closed subintervals of $I$ with $\lambda(U)\le2^{-k}$. For a fixed $U$, the length $\lambda(h(U))$ is a finite sum of absolute differences of values of $h$ at endpoints of intervals, and hence depends continuously on $h\in\mathcal{O}_J$. Thus $\mathcal{B}_{J,k}$ is closed in $\mathcal{O}_J$, and hence in $\Homeo(X)$.

We show that $\mathcal{B}_{J,k}$ has empty interior. Let $h\in\mathcal{O}_J$ and $\eta>0$. Divide $J$ into $M$ closed intervals $J_i$ of equal length, so fine that $|J|/M<\eta/2$ and
\[
|y-y'|\le|J|/M\implies|h^{-1}(y)-h^{-1}(y')|<\eta/2\qquad(y,y'\in J).
\]
In the interior of each $J_i$ choose a nondegenerate closed interval $W_i$ with $\lambda(h^{-1}(W_i))\le2^{-k}/M$; this is possible by the continuity of $h^{-1}$. Let $\psi\colon J\to J$ be the increasing piecewise affine homeomorphism that fixes the endpoints of every $J_i$ and maps each $W_i$ onto the closed central subinterval of $J_i$ of length $\frac34|J_i|$, affinely on $W_i$ and on the two remaining subintervals of $J_i$. Extend $\psi$ to $X$ by the identity outside $J$; since $\psi$ fixes the endpoints of $J$, this extension is a homeomorphism of $X$ by finite closed pasting.

Put $g=\psi\circ h$. Both $\psi$ and $\psi^{-1}$ map every $J_i$ onto itself, so the two mesh conditions give $\rho(g,h)<\eta$. For $U=\bigcup_ih^{-1}(W_i)$ we have
\[
\lambda(U)\le2^{-k},\qquad\lambda(g(U))=\tfrac34|J|,
\]
so $g\notin\mathcal{B}_{J,k}$. This proves that each $\mathcal{B}_{J,k}$ is nowhere dense. Consequently, the set
\[
\mathcal{G}_I=\Homeo(X)\setminus\bigcup_{J\in\mathcal{C}(X)}\bigcup_{k\ge1}\mathcal{B}_{J,k}
\]
is comeager.

Let $h\in\mathcal{G}_I$ and $J=h(I)$. For every $k$ choose a finite union $U_k\subset I$ of closed intervals with $\lambda(U_k)\le2^{-k}$ and $\lambda(h(U_k))>|J|/2$, and set
\[
N=\bigcap_{m\ge1}\bigcup_{k\ge m}U_k.
\]
Since $\sum_k\lambda(U_k)<\infty$, the Borel set $N$ is null. As $h$ is a bijection,
\[
h(N)=\bigcap_{m\ge1}\bigcup_{k\ge m}h(U_k).
\]
This is a decreasing intersection of measurable subsets of $J$, each of measure greater than $|J|/2$, so continuity of measure from above gives $\lambda(h(N))\ge|J|/2$. In particular, $h\notin\MC(X)$. Finally, the intersection of the sets $\mathcal{G}_I$ over the countable family $\mathcal{C}(X)$ is comeager and disjoint from $\MC(X)$.
\end{proof}

\subsection{Approximation and the support stabilizer}

\begin{lemma}[Relative approximation in the pure classes]\label{lem:approximation}
Let $A,B$ be Cantorvals whose boundaries are both null or both everywhere fat. Let $h\colon A\to B$ be a homeomorphism, and let $p_1,\dots,p_s\in A$ be distinct exposed points such that $h(p_1),\dots,h(p_s)$ are exposed in $B$. For every $\eta>0$ there is a measure-class homeomorphism $g\colon A\to B$ such that
\[
g(p_i)=h(p_i)\quad(1\le i\le s),\qquad\sup_A|g-h|+\sup_B|g^{-1}-h^{-1}|<\eta.
\]
The list of prescribed points may be empty.
\end{lemma}

\begin{proof}
Select the finitely many interval components $I$ of $A$ with $|I|\ge\eta/8$ or $|h(I)|\ge\eta/8$. Remove their interiors from $A$ and the interiors of their images from $B$; this gives sets $A_0$ and $B_0$ with $h(A_0)=B_0$. By Lemma~\ref{lem:clopen-pieces}, $A_0$ and $B_0$ are Cantorvals with the same boundaries as $A$ and $B$. The prescribed points remain exposed in $A_0$, their images remain exposed in $B_0$, and the endpoints of the removed intervals are exposed as well.

On $A_0$ consider the compatible metric
\[
d_h(u,v)=\max\{|u-v|,|h(u)-h(v)|\}.
\]
Every component of $A_0$ has $d_h$-diameter less than $\eta/8$. In a compact Hausdorff space, components coincide with quasicomponents. Hence every component has a clopen neighborhood of $d_h$-diameter less than $\eta/4$: take an open neighborhood $V$ of that diameter; since the component is the intersection of the clopen sets containing it, compactness of $A_0\setminus V$ yields finitely many such clopen sets whose intersection lies in $V$. A finite cover of $A_0$ by such neighborhoods, refined to a disjoint family, gives a clopen partition
\[
A_0=\bigsqcup_{i=1}^qC_i,\qquad\diam(C_i)<\eta/4,\qquad\diam(h(C_i))<\eta/4.
\]
By Lemma~\ref{lem:clopen-pieces}, every nonempty $C_i$ and every $h(C_i)$ is a Cantorval, whose boundary is of the same type as that of $A$ or $B$.

In each $C_i$, mark the prescribed points that it contains and the endpoints of removed intervals that belong to it; their images under $h$ are exposed in $h(C_i)$. Lemma~\ref{lem:matching} gives a measure-class homeomorphism $g_i\colon C_i\to h(C_i)$ that agrees with $h$ at all marked points. On each removed interval $I$, choose a partition so fine that its subintervals and their images under $h$ have lengths less than $\eta/4$. The piecewise affine interpolation of $h$ at the partition points is a monotone homeomorphism $g_I\colon I\to h(I)$ with the same endpoint values as $h$, and both $|g_I-h|$ and $|g_I^{-1}-h^{-1}|$ are less than $\eta/4$.

The finitely many maps $g_i$ and $g_I$ agree at common endpoints, and they paste to a bijection $g\colon A\to B$, which is a measure-class homeomorphism by finite closed pasting. On $C_i$, both $g$ and $h$ take values in $h(C_i)$, and for points of $h(C_i)$ both inverse maps take values in $C_i$. The diameter bounds for the partition and the interpolation bounds therefore give the stated estimate, and all marked values are preserved.
\end{proof}

\begin{lemma}\label{lem:orbit}
Let $X$ be a Cantorval. Any two points of $X$ whose connected components in $X$ are singletons lie in the same orbit of $\Homeo(X)$.
\end{lemma}

\begin{proof}
Let $p$ and $q$ be such points. Since components coincide with quasicomponents in compact spaces, $p$ and $q$ have neighborhood bases consisting of clopen sets. Choose strictly decreasing sequences of clopen sets $U_n\ni p$ and $V_n\ni q$, starting with $U_0=V_0=X$, whose diameters tend to $0$, so that $\bigcap_nU_n=\{p\}$ and $\bigcap_nV_n=\{q\}$. All the differences $U_n\setminus U_{n+1}$ and $V_n\setminus V_{n+1}$ are then nonempty and clopen, hence Cantorvals by Lemma~\ref{lem:clopen-pieces}, and therefore homeomorphic to each other (Section~\ref{sec:mc}). Choose homeomorphisms $U_n\setminus U_{n+1}\to V_n\setminus V_{n+1}$ for all $n$, paste them, and send $p$ to $q$. The resulting bijection of $X$ maps $U_n$ onto $V_n$ for every $n$. This gives continuity at $p$ and continuity of the inverse at $q$; continuity elsewhere holds because the differences $U_n\setminus U_{n+1}$ and $V_n\setminus V_{n+1}$ are open.
\end{proof}

\begin{theorem}\label{thm:density}
Let $X$ be a Cantorval, $K=\partial X$ and $S=S(X)$. Then
\[
\overline{\MC(X)}=\Homeo(X)\iff S=\varnothing\ \text{or}\ S=K.
\]
If $S$ is relatively clopen in $K$, then
\[
\overline{\MC(X)}=\{h\in\Homeo(X):h(S)=S\}.
\]
\end{theorem}

\begin{proof}
By Lemmas~\ref{lem:supp} and~\ref{lem:homeo-invariance}, every element of $\MC(X)$ maps $S$ onto itself. The setwise stabilizer of $S$ is closed in $\Homeo(X)$, because uniform convergence of maps and of their inverses preserves the two inclusions $h(S)\subset S$ and $h^{-1}(S)\subset S$. Hence $\overline{\MC(X)}$ is contained in this stabilizer.

Suppose that $\varnothing\ne S\ne K$. Then $S$ has positive measure, while $K\setminus S$ is a nonempty relatively open subset of a Cantor set; so both sets are uncountable. The set of endpoints of interval components is countable, so we can choose $p\in S$ and $q\in K\setminus S$ that are not such endpoints. Their components in $X$ are singletons, and by Lemma~\ref{lem:orbit} some homeomorphism of $X$ maps $p$ to $q$. Thus the stabilizer of $S$ is a proper closed subgroup, and $\MC(X)$ is not dense. If $S=\varnothing$ or $S=K$, Lemma~\ref{lem:approximation} with $A=B=X$ and no prescribed points shows that $\MC(X)$ is dense.

It remains to prove the reverse inclusion in the second formula when $S$ and $K\setminus S$ are both nonempty and clopen in $K$. Use the decomposition of Lemma~\ref{lem:bridges},
\[
X=X_S\cup X_R\cup B_1\cup\dots\cup B_m,
\]
into a fat part, a null part and finitely many bridges. A homeomorphism $h$ with $h(S)=S$ also preserves $K\setminus S$; since it preserves interval components and their endpoints, it maps each of the two parts onto itself and permutes the bridges. Given $\eta>0$, apply Lemma~\ref{lem:approximation} to the restrictions of $h$ to $X_S$ and to $X_R$, prescribing the values of $h$ at all endpoints of bridges (these points and their images are exposed in the two parts), with errors less than $\eta/2$. On each bridge, approximate $h$ by a monotone piecewise affine homeomorphism with the same endpoint values and with both uniform errors less than $\eta/2$. Finite closed pasting gives a measure-class homeomorphism $g$ of $X$ with $\rho(g,h)<\eta$. This proves the reverse inclusion; the pure cases were treated above.
\end{proof}

\begin{corollary}\label{cor:achievable-groups}
For every achievable Cantorval $E$, the group $\MC(E)$ is a dense meager subgroup of $\Homeo(E)$.
\end{corollary}

\begin{proof}
By the support dichotomy, Theorem~\ref{thm:dichotomy}, either $S(E)=\varnothing$ or $S(E)=\partial E$. Apply Theorems~\ref{thm:meager} and~\ref{thm:density}.
\end{proof}

\appendix

\section{Maps on co-small subspaces}\label{app:cosmall}

Let $(X,\mathcal{T},\mathcal{I})$ be a measure-class space. The members of $\mathcal{I}$ are called \emph{small}, and a subset $U\subset X$ is \emph{co-small} if $X\setminus U$ is small. All subspaces below carry the induced topology and the inherited ideal.

\begin{definition}\label{def:partial}
A \emph{partial measure-class homeomorphism} from $X$ to $Y$ is a homeomorphism $h\colon U\to V$ between co-small subspaces $U\subset X$ and $V\subset Y$ such that
\[
A\in\mathcal{I}\iff h(A)\in\mathcal{J}\qquad(A\subset U).
\]
Two such maps are identified if they agree on a co-small subset of $X$ contained in the intersection of their domains.
\end{definition}

\begin{proposition}\label{prop:groupoid}
Partial measure-class homeomorphisms, modulo co-small agreement, form a groupoid under restriction and composition.
\end{proposition}

\begin{proof}
Let $h\colon U\to V$ and $g\colon P\to Q$ be partial measure-class homeomorphisms from $X$ to $Y$ and from $Y$ to $Z$, respectively. Their composite is defined on
\[
W=h^{-1}(V\cap P),\qquad(g\circ h)(W)=g(V\cap P).
\]
The set $U\setminus W=h^{-1}(V\setminus P)$ is small in $X$, because $V\setminus P$ is small in $Y$; so $W$ is co-small. Likewise, $Q\setminus g(V\cap P)=g(P\setminus V)$ is small in $Z$. The composite is a homeomorphism between these co-small subspaces and preserves small sets in both directions.

If $h$ is replaced by a map that agrees with it on a co-small set, restrict both maps to that set. If $g$ is replaced by a map that agrees with it on a co-small subset of $Y$, discard the small complement of that subset and its preimage under $h$. Thus composition does not depend on the choice of representatives. Co-small agreement is an equivalence relation, since finite intersections of co-small sets are co-small. The image under $h$ of a co-small subset of $U$ is co-small in $Y$, which shows that inversion does not depend on the representatives either. Identities and associativity follow by restriction to common co-small domains.
\end{proof}

\begin{definition}\label{def:I-homotopy}
Continuous maps $f,g\colon X\to Y$ are \emph{$\mathcal{I}$-homotopic} if there are a co-small set $U\subset X$ and a continuous map $H\colon U\times[0,1]\to Y$ with $H(x,0)=f(x)$ and $H(x,1)=g(x)$ for $x\in U$. For metric spaces $X$ and $Y$, a map $f\colon X\to Y$ is \emph{$\mathcal{I}$-Lipschitz} with constant $L$ if it is $L$-Lipschitz on a co-small subspace. Self-maps $f,g\colon X\to X$ \emph{$\mathcal{I}$-commute} if $f\circ g=g\circ f$ on a co-small subspace.
\end{definition}

The relation of $\mathcal{I}$-homotopy is an equivalence relation: two homotopies can be restricted to a common co-small subspace and concatenated there.

\begin{lemma}\label{lem:dense-cosmall}
Assume that every co-small subset of $X$ is dense. If $Y$ is Hausdorff, then continuous maps $f,g\colon X\to Y$ that agree on a co-small set agree everywhere. If $X$ and $Y$ are metric spaces and a continuous map $f\colon X\to Y$ is $\mathcal{I}$-Lipschitz with constant $L$, then it is $L$-Lipschitz on $X$. In particular, continuous self-maps of a Hausdorff space that $\mathcal{I}$-commute commute everywhere.
\end{lemma}

\begin{proof}
The set where two continuous maps into a Hausdorff space agree is closed, and here it contains a dense set. For the metric assertion, let $U$ be a co-small subspace on which $f$ is $L$-Lipschitz; then $U\times U$ is dense in $X\times X$, and the Lipschitz inequality passes to the limit.
\end{proof}

The density assumption holds for the meager ideal of a Baire space and for the null ideal of a measure of full support. The following construction works for both ideals at once in Euclidean space.

\begin{proposition}\label{prop:euclidean}
For every $n\ge1$ there are co-null, co-meager subspaces of $\R^n$ and of $\R$ and a homeomorphism between them that preserves both the inherited Lebesgue null ideals and the inherited meager ideals.
\end{proposition}

\begin{proof}
Let $\Omega=\{0,1\}^{\mathbb{N}}$ carry the product topology and the fair Bernoulli probability measure $\nu$. Removing the eventually constant sequences, we obtain a co-null, co-meager subspace $B$. The binary expansion map
\[
\pi((s_k)_{k\ge1})=\sum_{k\ge1}s_k2^{-k}
\]
restricts to a homeomorphism from $B$ onto the set $P$ obtained from $(0,1)$ by removing the dyadic rationals. Indeed, a finite binary prefix determines the intersection of an open dyadic interval with $P$, and these sets show that the inverse map is continuous. The map $\pi$ sends $\nu$ to Lebesgue measure; this is checked first on dyadic intervals and then extends to all Borel sets.

Interleaving the coordinates defines a measure-preserving homeomorphism $T_n\colon\Omega^n\to\Omega$,
\[
\bigl(T_n(s^{(1)},\dots,s^{(n)})\bigr)_{n(k-1)+i}=s^{(i)}_k\qquad(k\ge1,\ 1\le i\le n).
\]
Its image $C_n=T_n(B^n)$ is co-null and co-meager in $\Omega$, and it is contained in $B$. Consequently,
\[
\pi\circ T_n\circ(\pi^{-1})^n\colon P^n\longrightarrow\pi(C_n)
\]
is a homeomorphism between co-null, co-meager subspaces of $(0,1)^n$ and $(0,1)$. It preserves the Lebesgue measure of Borel sets, and taking Borel null hulls shows that it preserves null sets in both directions. The inherited meager ideals coincide with the meager ideals of these dense subspaces, and are therefore preserved by the homeomorphism.

Finally, conjugate by a smooth diffeomorphism $\R\to(0,1)$ and by its coordinatewise product. These maps and their inverses are locally Lipschitz and hence preserve null sets; being homeomorphisms, they also preserve meager sets.
\end{proof}

\begin{example}\label{ex:circle}
On $S^1$, with either normalized arc length or the meager ideal, the maps $z\mapsto z^k$ and $z\mapsto z^\ell$ are homotopic on the complement of one point, for all $k,\ell\in\mathbb{Z}$. Indeed, on $S^1\setminus\{1\}$ write $z=e^{i\theta}$ with $0<\theta<2\pi$ and put
\[
H(e^{i\theta},t)=\exp\bigl(i((1-t)k+t\ell)\theta\bigr).
\]
This is a continuous homotopy on the co-small subspace $S^1\setminus\{1\}$. On the whole circle, of course, the ordinary homotopy classes are still distinguished by their degrees.
\end{example}

\section{Generic affine-line chain homotopy}\label{app:generic}

Let $k$ be an infinite field, let $X$ be an integral separated scheme of finite type over $k$, and let $Y$ be a separated scheme of finite type over $k$. Write $K=k(X)$ and $Y_K=Y\times_k\Spec K$. For a $k$-algebra $B$, put $Y(B)=\Hom_k(\Spec B,Y)$.

\begin{definition}\label{def:generic-homotopy}
A subset of $X$ is \emph{thin} if it is contained in a closed subset whose complement is dense. Two morphisms $f,g\colon U\to Y$, where $U\subseteq X$ is a dense open subset, are \emph{elementarily $\mathbb{A}^1$-homotopic} if there is a morphism
\[
H\colon U\times_k\mathbb{A}^1_k\longrightarrow Y\qquad\text{with}\qquad H(-,0)=f,\quad H(-,1)=g.
\]
Two rational maps $\phi,\psi\colon X\dashrightarrow Y$ are \emph{generically $\mathbb{A}^1$-chain homotopic} if, on some common dense open subset $U$, their representatives are joined by a finite chain of such homotopies. The set of equivalence classes is denoted by $\Pi(X,Y)$.
\end{definition}

The thin subsets form a proper ideal, that is, a proper family closed under taking subsets and finite unions. This ideal need not be closed under countable unions, since a countably infinite set of $k$-points of $\mathbb{A}^1_k$ is dense. Reversing the parameter of a homotopy gives symmetry, and intersecting the finitely many domains of two chains gives transitivity of the relation in Definition~\ref{def:generic-homotopy}.

For a field extension $L/k$, let $\sim_L$ be the equivalence relation on $Y(L)$ generated by the pairs $(\gamma(0),\gamma(1))$, where $\gamma\colon\mathbb{A}^1_L\to Y_L$ is an $L$-morphism. This is the pointwise affine-line chain relation used in the study of $\mathbb{A}^1$-chain connectedness; see \cite[Section~2.2]{AM}.

\begin{proposition}\label{prop:generic-point}
Restriction to the generic point induces a natural bijection
\[
\Pi(X,Y)\simeq Y(K)/{\sim_K}.
\]
In particular, for fixed $Y$ the set $\Pi(X,Y)$ depends only on the function field $k(X)$ over $k$.
\end{proposition}

\begin{proof}
Choose a nonempty affine open subset $X_0=\Spec A$ of $X$. The principal open sets $U_f=\Spec A_f$, $0\ne f\in A$, are cofinal among the dense open subsets of $X$ after shrinking, and
\[
K=\varinjlim_{f\ne0}A_f,\qquad K[s]=\varinjlim_{f\ne0}A_f[s].
\]
Since $Y$ is of finite presentation over $k$, the limit theorem for morphisms of finite presentation \cite[Th{\'e}or{\`e}me~8.8.2]{EGA} gives
\begin{align*}
Y(K)&\simeq\varinjlim_{f\ne0}Y(A_f),\\
\Hom_K(\mathbb{A}^1_K,Y_K)&\simeq\varinjlim_{f\ne0}\Hom_k(U_f\times_k\mathbb{A}^1_k,Y).
\end{align*}
Thus every $K$-point extends to a rational map on $X$, and every affine-line path over $K$ extends to a homotopy on some $U_f$. The first bijection identifies rational maps with their generic values; the separatedness of $Y$ also shows directly that representatives with the same generic value agree on their common domain.

A generic homotopy chain restricts to a chain over $K$. Conversely, extend the finitely many points and paths of a chain over $K$, including the given representatives of its endpoints, and then shrink to one common $U_f$. The endpoint identities hold after a further shrinking, by the same limit theorem. The resulting paths form the required chain on $U_f$, which proves both directions of the identification.
\end{proof}

\begin{proposition}\label{prop:R-equivalence}
If $Y$ is proper over $k$, then $\sim_K$ is Manin's $R$-equivalence, and
\[
\Pi(X,Y)\simeq Y(K)/R.
\]
\end{proposition}

\begin{proof}
Direct $R$-equivalence relates the values at $0$ and $1$ of a rational map $\mathbb{P}^1_K\dashrightarrow Y_K$ that is defined at these points, and $R$-equivalence is the equivalence relation generated by it~\cite{Manin}. Such a rational map extends across each missing point of $\mathbb{P}^1_K$: the local ring at such a point is a discrete valuation ring, and the valuative criterion of properness provides the extension. The extensions glue by separatedness, and restricting the resulting morphism to $\mathbb{A}^1_K$ gives an affine-line path. Conversely, every affine-line path defines a rational map on $\mathbb{P}^1_K$ that is defined at $0$ and $1$. Hence the generating relations coincide, and so do the equivalence relations they generate.
\end{proof}

\begin{example}[Projective space]\label{ex:projective}
For every $n\ge0$, the set $\Pi(X,\mathbb{P}^n_k)$ has exactly one element. Indeed, any two $K$-points of projective space lie in a common affine chart, since $K$ is infinite and hence there is a linear form that vanishes at neither point. In that chart the affine line $s\mapsto(1-s)p+sq$ joins them. In particular, for every integer $d\ge1$ the formula
\[
H(t,s)=(1-s)t^d+st\qquad\text{on}\quad\mathbb{A}^1_t\times_k\mathbb{A}^1_s
\]
gives a generic homotopy, with target $\mathbb{A}^1_k\subset\mathbb{P}^1_k$, from $t\mapsto t^d$ to the identity. By contrast, the degree is an invariant in the classification of pointed naive homotopy classes of maps defined on all of $\mathbb{P}^1_k$; see \cite[Proposition~2.3]{Cazanave}.
\end{example}

\begin{example}[Four arithmetic classes]\label{ex:elliptic}
Let $C/\mathbb{Q}$ be the smooth projective elliptic curve with affine equation $y^2=x^3-x$ and point at infinity $O$, and put $Y=C\times_{\mathbb{Q}}\mathbb{P}^1_{\mathbb{Q}}$. Then
\[
\Pi(\mathbb{P}^1_{\mathbb{Q}},Y)\simeq C(\mathbb{Q}(t))=C(\mathbb{Q})=\{O,(0,0),(1,0),(-1,0)\}.
\]
To prove the first identification, set $L=\mathbb{Q}(t)$. A morphism $\mathbb{A}^1_L\to C_L$ extends to $\mathbb{P}^1_L$ by properness and is constant: after base change to an algebraic closure, a nonconstant morphism from a curve of genus zero to a curve of genus one is impossible by the Riemann--Hurwitz formula in characteristic zero. Thus a chain in $Y(L)$ preserves the $C$-coordinate, while each fiber over a point of $C(L)$ is chain connected by Example~\ref{ex:projective}. Proposition~\ref{prop:generic-point} now applies.

An element of $C(\mathbb{Q}(t))$ similarly gives a rational map $\mathbb{P}^1_{\mathbb{Q}}\dashrightarrow C$, which extends to a morphism and is constant; evaluation at a rational point shows that the constant lies in $C(\mathbb{Q})$. Finally, if $(x,y)\in C(\mathbb{Q})$ and $y\ne0$, then
\[
\Bigl|\frac{x^2-1}{y}\Bigr|,\qquad\Bigl|\frac{2x}{y}\Bigr|,\qquad\Bigl|\frac{x^2+1}{y}\Bigr|
\]
are the sides of a rational right triangle of area $|x(x^2-1)|/y^2=1$. By Fermat's right-triangle theorem, $1$ is not a congruent number \cite[Chapter~I, p.~2]{Koblitz}. Hence $y=0$ and $x\in\{-1,0,1\}$; together with $O$ these are all the rational points.
\end{example}

\subsection*{Acknowledgments}
The author is grateful to Piotr Nowakowski for reading the manuscript carefully, for many helpful comments and suggestions, and for pointing out several references.

\end{document}